\documentclass[11pt]{amsart}

\usepackage{amsmath}
\usepackage{amssymb}
\usepackage{comment}
\usepackage{color}
\usepackage[colorlinks,citecolor=blue,urlcolor=black,linkcolor=black]{hyperref}

\newtheorem{theorem}{Theorem}[section]
\newtheorem{claim}[theorem]{Claim}

\newtheorem{lemma}[theorem]{Lemma}

\theoremstyle{definition}
\newtheorem{definition}[theorem]{Definition}

\newtheorem{question}[theorem]{Question}

\theoremstyle{remark}

\newcount\skewfactor
\def\mathunderaccent#1#2 {\let\theaccent#1\skewfactor#2
\mathpalette\putaccentunder}
\def\putaccentunder#1#2{\oalign{$#1#2$\crcr\hidewidth
\vbox to.2ex{\hbox{$#1\skew\skewfactor\theaccent{}$}\vss}\hidewidth}}

\def\smallbox#1{\leavevmode\thinspace\hbox{\vrule\vtop{\vbox
   {\hrule\kern1pt\hbox{\vphantom{\tt/}\thinspace{\tt#1}\thinspace}}
   \kern1pt\hrule}\vrule}\thinspace}

\newcommand{\cf}{{\rm cf}}

\def\qedref#1{$\qed_{\reforiginal{#1}}$}

\title{Tiltan and graphs with no infinite paths}
\author{Shimon Garti}
\address{Institute of Mathematics,
 The Hebrew University of Jerusalem,
 Jerusalem 91904, Israel}
\email{shimon.garty@mail.huji.ac.il}
\thanks{}


\begin{document}
\let\labeloriginal\label
\let\reforiginal\ref

\begin{abstract}
We prove the consistency of tiltan with the positive relation $\omega^*\cdot\omega_1 \rightarrow (\omega^*\cdot\omega_1,\text{infinite path})^2$.
\end{abstract}

\maketitle

\newpage

\section{Introduction}

Let $G=(V,E)$ be a graph.
An independent subset of $V$ is a set of vertices $W\subseteq V$ such that $[W]^2\cap E=\varnothing$.
An infinite path in $G$ is a sequence of vertices $\langle v_n:n\in\omega \rangle$ with no repetitions such that $\{v_n,v_{n+1}\}\in E$ for every $n\in\omega$.
Intuitively, these concepts are orthogonal.
If one wishes to eliminate large independent sets then one must add edges to many pairs. In such cases it becomes harder to avoid infinite paths.
For making this intuition precise we need a definition of \emph{large} independent sets.
The most natural suggestion would be a subset $W$ of $V$ with the same order type.

\begin{definition}
\label{defarrow} The relation $\tau\rightarrow (\tau,\text{infinite path})^2$ means that for every graph $G=(V,E)$ with ${\rm otp}(V)=\tau$ there exists either an independent set $W\subseteq V$ so that ${\rm otp}(W)=\tau$ or an infinite path.
\end{definition}

By order type we do not confine ourselves to well-orderings. Rather, we refer to a variety of structures. We consider ordinals $\alpha$ with their well-ordering, the backward ordering $\alpha^*$ and ordinal products of these types.
All graphs in this paper are undirected.

The notation $\tau\rightarrow (\tau,\text{infinite path})^2$ comes from partition theorems of infinite combinatorics.
Given a graph $G$ one may think of a coloring of its pairs with two colors.
The first color is assigned to every pair of vertices with no edge, and the second color is given to pairs with an edge.
The positive relation states that there exists a full sized subset with the first color or an infinite sequence with the second color.

We shall focus on the order type $\omega^*\cdot\omega_1$.
For a convenient and concrete example, if the ambient set is $\omega_1\times\omega$ then the order defined by $(\alpha,m)<^*(\beta,n)$ iff $(\alpha<\beta)\vee(\alpha=\beta$ and $m>n)$ is of type $\omega^*\cdot\omega_1$.
A convenient way to visualize this type is by thinking about $\omega_1$ many columns, each of which is a copy of $\omega^*$.
Regarding the above definition one may wonder whether $\omega^*\cdot\omega_1 \rightarrow (\omega^*\cdot\omega_1,\text{infinite path})^2$.

In the parallel abstract situation of infinite combinatorics, $\lambda\rightarrow(\lambda,\omega)^2$ for every infinite cardinal $\lambda$, and this is known as the Erd\H{o}s-Dushnik-Miller theorem.
However, if $\alpha$ is an ordinal but not a cardinal then $\alpha\nrightarrow(\alpha,\omega)^2$.
These facts motivate the investigation of more types like $\omega^*\cdot\omega_1$.
We indicate that an infinite path in a graph is a weaker notion than an infinite monochromatic set, since the homogeneity is required only at consecutive elements of the path.
There is some evidence that the existence of monochromatic paths is strictly weaker than the existence of monochromatic sets, see \cite{MR4245063} and \cite{MR4217964}.
In our context, one may obtain such paths even though the order-type of the graph is neither a cardinal, nor an ordinal.
Put another way, a mysterious path may show up, as described in \cite[page 10]{comet}: Moomintroll was just putting up a swing when Sniff got home. He seemed very interested in the mysterious path, and directly after lunch they set off to have a look at it.

Back to the context of graph theory, the above relation cannot be decided by the axioms of set theory.
Namely, one can prove the consistency of $\omega^*\cdot\omega_1 \rightarrow (\omega^*\cdot\omega_1,\text{infinite path})^2$ in some extension of ZFC on the one hand, and one can show that $\omega^*\cdot\omega_1 \nrightarrow (\omega^*\cdot\omega_1,\text{infinite path})^2$ in another extension on the other hand.
The negative direction was given by Baumgartner and Larson, \cite{MR1050558}, and the positive direction by Larson in \cite{MR1050560}.
It is done in these papers through the classical way of confronting the constructible universe with the universe under Martin's axiom with large continuum.

Actually, the full strength of the constructible universe is not required.
Baumgartner and Larson constructed a graph $G=(V,E)$ of type $\omega^*\cdot\omega_1$ with no independent subset of this type and no infinite path merely from the diamond principle at $\aleph_1$.
Recall that $\Diamond_{\aleph_1}$ says that there exists a sequence of sets $\langle A_\alpha:\alpha\in\omega_1\rangle$ such that $A_\alpha\subseteq\alpha$ for every $\alpha\in\omega_1$ and for every $A\subseteq\omega_1$ the set $S_A = \{\alpha\in\omega_1:A\cap\alpha=A_\alpha\}$ is a stationary subset of $\omega_1$.
The opposite direction employs Martin's axiom with $2^\omega>\omega_1$, and then $\omega^*\cdot\omega_1 \rightarrow (\omega^*\cdot\omega_1,\text{infinite path})^2$. Both directions are elaborated in another paper of Larson, \cite{MR2279658}.
In this paper she explains the importance of the type $\omega^*\cdot\omega_1$ and poses the problem which stands in the hub of this paper.
We let the following definition into the discussion.

\begin{definition}
\label{deftiltan} Tiltan. \newline
Let $\kappa=\cf(\kappa)>\aleph_0$. \newline
The tiltan principle $\clubsuit_\kappa$ says that there exists a sequence $\langle T_\alpha:\alpha$ is a limit ordinal of $\kappa\rangle$ such that each $T_\alpha$ is a cofinal subset of $\alpha$ and for every $A\in[\kappa]^\kappa$ the set $S_A = \{\alpha\in\kappa:T_\alpha\subseteq A\}$ is a stationary subset of $\kappa$.
\end{definition}

The common name of this statement is the club principle.
It has been introduced by Ostaszewski, in \cite{MR0438292}.
We shall call it tiltan\footnote{Let us indicate that in some good old manuscripts the pronunciation is \emph{taltan}, see the relevant discussion in \cite{Yeivin}.} since the word \emph{club} is extensively used as an acronym for closed and unbounded sets.

The tiltan follows from the diamond, and it is strictly weaker than the diamond.
In particular, $\Diamond_{\aleph_1}\Rightarrow 2^\omega=\omega_1$ while $\clubsuit_{\aleph_1}$ is consistent with $2^\omega>\omega_1$.
Remark that Martin's axiom with $2^\omega>\omega_1$ implies $\neg\clubsuit_{\aleph_1}$.
Therefore, the following question of Larson from \cite{MR2279658} is natural:

\begin{question}
\label{qlarson} Is it consistent that tiltan holds at $\aleph_1$ and concomitantly $\omega^*\cdot\omega_1 \rightarrow (\omega^*\cdot\omega_1,\text{infinite path})^2$?
\end{question}

We shall give a positive answer to this question.
Let us indicate that having a positive result is a bit surprising.
One of the main differences between tiltan and diamond is that the diamond prediction is based on equality ($A\cap\alpha=A_\alpha$) while the tiltan prediction gives only inclusion ($T_\alpha\subseteq A$).
In the negative arrow relation proved under diamond in \cite{MR1050558}, only inclusion is needed for the construction of a graph exemplifying $\omega^*\cdot\omega_1 \nrightarrow (\omega^*\cdot\omega_1,\text{infinite path})^2$.
Despite this fact, tiltan is consistent with the positive relation $\omega^*\cdot\omega_1 \rightarrow (\omega^*\cdot\omega_1,\text{infinite path})^2$, as we shall see.

The rest of the paper is arranged in two additional sections.
In the first one we unfold some background material and we try to explicate the idea behind the proof.
In the second, we prove the main theorem.
Our notation is mostly standard.
Let us mention the notation $S^\lambda_\kappa$ which refers to the set $\{\delta<\lambda:\cf(\delta)=\kappa\}$ where $\kappa$ is a regular cardinal.
We employ the Jerusalem forcing notation, so $p\leq q$ reads $p$ is weaker than $q$.
Consequently we shall speak about a least upper bound of conditions, a downward closed generic set, and so forth.
If $p$ is compatible with $q$ then we write $p\parallel q$.
If $p$ and $q$ are incompatible then we shall write $p\perp q$.
The meaning of the symbol $\exists^\infty$ is that there are infinitely many elements which satisfy the statement which falls under the scope of this quantifier.
We employ this notation with respect to sets of natural numbers.

I am deeply indebted to the anonymous referee for many mathematical corrections and a lot of helpful suggestions.
The referee pointed out a major flaw in the original version of the manuscript and enabled me to fix the problematic issue.
I learned several mathematical things from the work of the referee on my paper, but I learned much more from his/her infinite patience for paths and infinite path of patience.

\newpage

\section{Background}

Larson proved the consistency of $\omega^*\cdot\omega_1 \rightarrow (\omega^*\cdot\omega_1,\text{infinite path})^2$ under Martin's axiom.
A central component in our proof is a similar theorem at the level of $\omega_2$.
We shall use a generalized form of Martin's axiom, and there are several such theorems in the literature.
The most appropriate among which for our proof is Shelah's version from \cite{MR0505492}.

\begin{theorem}
\label{thmgma} Generalized Martin's axiom. \newline
One can force $2^{\aleph_0}=\aleph_1, 2^{\aleph_1}>\aleph_2$ and if $\mathbb{P}$ satisfies:
\begin{enumerate}
\item [$(a)$] If $p,q\in\mathbb{P}$ and $p\parallel q$ then there is a least upper bound for $p,q$ in $\mathbb{P}$.
\item [$(b)$] If $\langle p_i:i\in\omega\rangle$ is an increasing sequence of conditions in $\mathbb{P}$ then it has a least upper bound in $\mathbb{P}$.
\item [$(c)$] If $\{p_\alpha:\alpha\in\omega_2\}\subseteq\mathbb{P}$ then there is a club $C\subseteq\omega_2$ and a regressive function $f:\omega_2\rightarrow\omega_2$ such that for every $\alpha,\beta\in C\cap S^{\aleph_2}_{\aleph_1}$ if $f(\alpha)=f(\beta)$ then $p_\alpha\parallel p_\beta$.
\end{enumerate}
then for every $\kappa<2^{\aleph_1}$ such that $\gamma<\kappa\Rightarrow \gamma^{\aleph_0}<\kappa$ and any collection $\mathcal{D} = \{D_\eta: \eta\in\kappa\}$ of dense subsets of $\mathbb{P}$ there exists a filter $G\subseteq\mathbb{P}$ so that $G\cap D_\eta\neq\varnothing$ for every $\eta\in\kappa$.
\end{theorem}

\hfill \qedref{thmgma}

The forcing conditions in Larson's proof are finite independent sets.
In our proof the conditions are countable.
As a first step we shall use the generalized Martin's axiom in order to force $\omega^*\cdot\omega_2 \rightarrow (\omega^*\cdot\omega_2,\text{infinite path})^2$, and requirement $(b)$ above forces us to force with countable conditions.
This is one major difference between Martin's axiom and the generalized Martin's axiom which complicates the density argument.

Another problem is the chain condition.
Martin's axiom applies to any $ccc$ forcing notion, but all the generalizations to higher cardinals require more than $\kappa$-cc, and it is known that $\kappa$-cc is insufficient.
In Shelah's version, the strengthening of the chain condition is reflected in requirement $(c)$.
Our proof of $(c)$ in the specific forcing notion of this paper is based on the ordinary partition relation $\omega_2\rightarrow(\omega_2-st,\omega_1)^2$ which says that for every coloring $d:[\omega_2]^2\rightarrow\{0,1\}$ one can find either a $1$-monochromatic sequence of type $\omega_1$ or a stationary $0$-monochromatic subset of $\omega_2$.
The following is a folklore but we give the proof since we will use both the statement and the argument within the proof.

\begin{lemma}
\label{lemstat} Assume $2^\omega=\omega_1$.
Then $\omega_2\rightarrow(\omega_2-st,\omega_1)^2$.
Moreover, $T\rightarrow(\omega_2-st,\omega_1)^2$ whenever $T\subseteq S^{\omega_2}_{\omega_1}$ is stationary.
\end{lemma}

\par\noindent\emph{Proof}. \newline
Let $d:[\omega_2]^2\rightarrow\{0,1\}$ be a coloring.
If there is a $1$-monochromatic sequence of length $\omega_1$ then we are done.
Suppose that there is no such a sequence.
For every $\delta\in S^{\omega_2}_{\omega_1}$ choose a sequence $c_\delta$ of ordinals below $\delta$ such that $c_\delta\cup\{\delta\}$ is $1$-monochromatic and $c_\delta$ is maximal with this property.

By our assumption, $c_\delta$ is bounded below $\delta$ since $\cf(\delta)=\omega_1$.
Hence the mapping $h(\delta)=\sup(c_\delta)$ is regressive on $S^{\omega_2}_{\omega_1}$.
Choose $\eta\in\omega_1$ and a stationary $S\subseteq S^{\omega_2}_{\omega_1}$ such that $h(\delta)=\eta$ for every $\delta\in{S}$.
Notice that $\eta<\min(S)$.
Since $2^\omega=\omega_1$, there are only $\aleph_1$-many sequences of the form $c_\delta$ (recall that $\eta$ is an ordinal less than $\omega_2$ and each $c_\delta$ is a subset of $\eta$).
Hence by shrinking $S$ if needed we may assume that there is a fixed sequence $c$ such that $c_\delta=c$ for every $\delta\in{S}$.

We claim that $S$ is $0$-monochromatic under $d$.
To see this, suppose that $\delta,\varepsilon\in{S}$ and $\delta<\varepsilon$.
If $d(\delta,\varepsilon)=1$ then $c\cup\{\delta\}$ is $1$-monochromatic with $\varepsilon$ and then $h(\varepsilon)\geq\delta>\eta$.
This is impossible since $\varepsilon\in{S}$.
Hence necessarily $d(\delta,\varepsilon)=0$ whenever $\{\delta,\varepsilon\}\subseteq{S}$, so we are done.
The additional part of the lemma is proved in the same way, upon replacing $S^{\omega_2}_{\omega_1}$ by $T$.

\hfill \qedref{lemstat}

The next issue is a special kind of tiltan which we shall need for our proof.
Definition \ref{deftiltan} is phrased with respect to $\kappa$, but one can replace $\kappa$ by any stationary subset $S\subseteq\kappa$.
Clearly, if $S_0\subseteq S_1$ are stationary then $\clubsuit_{S_0} \Rightarrow \clubsuit_{S_1}$ and hence $\clubsuit_S \Rightarrow \clubsuit_\kappa$ whenever $S$ is a stationary subset of $\kappa$.
The following theorem from \cite{MR1623206} served for proving the consistency of tiltan at $\aleph_1$ with $2^\omega>\omega_1$, and we shall use it with respect to infinite graphs.

\begin{theorem}
\label{thmindest} Assume that $\Diamond_S$ holds at every stationary subset $S$ of $\aleph_1$ and $\aleph_2$.
Then one can define a tiltan sequence on $S^{\aleph_2}_{\aleph_0}$ which is indestructible upon any further forcing extension with an $\aleph_1$-complete forcing notion.
\end{theorem}

\hfill \qedref{thmindest}

We indicate that the proof of the generalized Martin's axiom employs an $\aleph_1$-complete forcing notion, hence preserves instances of indestructible tiltan. We shall use this fact in the proof of the main theorem.

We mention three additional classical theorems, to be used within our proof.
Firstly, Ramsey's theorem which says that $\omega\rightarrow(\omega)^2_\ell$ for every $\ell\in\omega$.
Namely, any coloring $c:[\omega]^2\rightarrow\ell$ admits a monochromatic infinite set.
Secondly, Hajnal's free set theorem which says that if $\kappa<\lambda, |A|=\lambda, f:A\rightarrow\mathcal{P}(A)$ is a set-mapping (i.e. $a\notin f(a)$ for every $a\in A$) and $|f(a)|<\kappa$ for every $a\in A$ then there exists an $f$-free subset $B\subseteq A$ of size $\lambda$. Recall that $B$ is $f$-free iff $B\cap f(b)=\varnothing$ whenever $b\in B$.
For the third theorem recall that if $\kappa$ is an infinite cardinal then ${\rm log}_\kappa(\kappa^+)=\min\{\theta:\kappa^\theta>\kappa\}$.
One can show that if $\kappa\geq\omega$ then $\kappa^+\rightarrow(\kappa^+,{\rm log}_\kappa(\kappa^+)+1)^2$, see \cite{MR795592}.
In particular, if $2^\omega=\omega_1$ then ${\rm log}_{\omega_1}(\omega_2)=\omega_1$ and hence $\omega_2\rightarrow(\omega_2,\omega_1)^2$.
In fact, one has the stronger relation $\omega_2\rightarrow(\omega_2-st,\omega_1)^2$ as proved above.

We shall also need a statement concerning path relations in the following weak form.
Call a coloring $d:\kappa\times\kappa\rightarrow\omega\times\omega$ \emph{anti-symmetric} iff $d(\alpha,\beta)=(i,j)\Leftrightarrow d(\beta,\alpha)=(j,i)$ whenever $\alpha,\beta\in\kappa$.
Let us say that $\kappa\rightarrow_{\rm asp}(\omega)^2_{\omega\times\omega}$ iff for every anti-symmetric coloring $d:\kappa\times\kappa\rightarrow\omega\times\omega$ one can find an infinite path $\psi=(\alpha_m:m\in\omega)$, the elements of $\psi$ being ordinals of $\kappa$ and for every $m\in\omega$ if $d(\alpha_m,\alpha_{m+1})=(i,j)\wedge d(\alpha_{m+1},\alpha_{m+2})=(k,\ell)$ then $j=k$.

\begin{lemma}
\label{lemwp} $\omega_1\rightarrow_{\rm asp}(\omega)^2_{\omega\times\omega}$.
\end{lemma}

\par\noindent\emph{Proof}. \newline
Let $d:\omega_1\times\omega_1\rightarrow\omega\times\omega$ be anti-symmetric.
Let $\chi$ be a sufficiently large regular cardinal and choose a countable elementary submodel $M\prec\mathcal{H}(\chi)$ so that $d\in M$.
Let $\delta=M\cap\omega_1$ be the characteristic ordinal of $M$ and notice that $\cf(\delta)=\omega$.

Fix an ordinal $\alpha_0\in\delta$ and assume that $d(\alpha_0,\delta)=(i,j)$.
Denote the set $\{\alpha\in\delta:d(\alpha,\delta)=(i,j)\}$ by $B$ and notice that $B$ is unbounded in $\delta$ by elementarity.
By definition, $\alpha_0\in B$.
Choose $\alpha_1>\alpha_0$ so that $\alpha_1\in B$.
This means that $d(\alpha_0,\delta)=d(\alpha_1,\delta)=(i,j)$ so by elementarity one can find $\beta_1>\alpha_1$ such that $\beta_1<\delta$ and $d(\alpha_0,\beta_1)=d(\alpha_1,\beta_1)=(i,j)$.
We choose now another element $\alpha_2\in B$ so that $\alpha_2>\beta_1$.
Since $d(\alpha_2,\delta)=(i,j)$ one can choose $\beta_2<\delta$ such that $\beta_2>\alpha_2$ and $d(\alpha_1,\beta_2)=d(\alpha_2,\beta_2)=(i,j)$.
We render this process in the same way by induction on $n\in\omega$ and finally define:
\[
\psi=(\alpha_0,\alpha_n,\beta_n:0<n<\omega).
\]
One can verify that $\psi$ forms an infinite path in the sense defined before the statement of the lemma.
We conclude, therefore, that $\omega_1\rightarrow_{\rm asp}(\omega)^2_{\omega\times\omega}$ as required.

\hfill \qedref{lemwp}

\newpage

\section{Graphs with no infinite path}

In this section we prove the main result of the paper, namely tiltan is consistent with $\omega^*\cdot\omega_1 \rightarrow (\omega^*\cdot\omega_1,\text{infinite path})^2$.
Let us describe the architecture of the proof.
The first step is to fix an indestructible tiltan sequence at $S=S^{\aleph_2}_{\aleph_0}$.
The second step is to force the generalized Martin's axiom so that $2^\omega=\omega_1, 2^{\omega_1}>\omega_2$ and the tiltan from the first step is preserved.
The main theorem at this stage is the positive relation $\omega^*\cdot\omega_2 \rightarrow (\omega^*\cdot\omega_2,\text{infinite path})^2$.
This relation will follow from the generalized Martin's axiom.
The final step is to collapse $\aleph_1$ by making it a countable ordinal.

It is easy to show that the tiltan is preserved by this collapse, in the sense that it holds in the generic extension over some stationary subset of $\aleph_1$.
Likewise, the above positive relation obtained by the generalized Martin's axiom becomes after the collapse $\omega^*\cdot\omega_1 \rightarrow (\omega^*\cdot\omega_1,\text{infinite path})^2$.
This general plan has been used by Shelah, \cite{MR1623206}, in his proof of the consistency of tiltan with $2^\omega>\omega_1$.
At the end of the paper we shall try to explain what are the features of a statement that one should expect to hold (consistently) with tiltan.

We commence with the concept of clean columns, as defined by Larson.
The definitions and claims in our context are adapted to the level of $\aleph_2$.
In the definition and lemma below we follow in the footsteps of Larson.

\begin{definition}
\label{defclean} Clean columns. \newline
Let $G=(V,E)$ be a graph where $V\subseteq\omega_2\times\omega$.
\begin{enumerate}
\item [$(\aleph)$] For every $\beta\in\omega_2$, the $\beta$th column of the graph is the set $V(\beta) = V\cap(\{\beta\}\times\omega)$.
\item [$(\beth)$] $G$ has clean columns iff the following three properties hold for every $\beta\in\omega_2$:
\begin{enumerate}
\item [$(a)$] Either $V(\beta)=\varnothing$ or $|V(\beta)|=\aleph_0$.
\item [$(b)$] $[V(\beta)]^2\cap E=\varnothing$.
\item [$(c)$] For every $(\alpha,n)\in V$ there is at most one pair $(\beta,\ell)$ such that $\{(\alpha,n),(\beta,\ell)\}\in E$.
\end{enumerate}
\end{enumerate}
\end{definition}

Graphs with clean columns simplify considerably the treatment of independent subsets and related notions. Of course, a graph $G$ may lack this property.
However, we focus on graphs of type $\omega^*\cdot\omega_2$ with no infinite path. In such graphs one can always pass to a subgraph of the same order type with clean columns.
Ahead of proving this assertion, we need a simple lemma.

\begin{lemma}
\label{lemfinite} Let $G=(V,E)$ be a graph over $\omega_2\times\omega$ with no infinite path, and let $\beta\in\omega_2$. Assume that $C\subseteq\{\beta\}\times\omega$ and $|C|=\aleph_0$. \newline
There exists a finite set $A\subseteq\omega_2\times\omega$ and an infinite set $B\subseteq C$ such that:
\begin{enumerate}
\item [$(a)$] Every element of $A$ is connected with every element of $B$.
\item [$(b)$] If $(\alpha,m)\notin A$ then there is at most one element $(\beta,n)\in B$ such that $\{(\alpha,m),(\beta,n)\}\in E$.
\item [$(c)$] $[B]^2\cap E=\varnothing$.
\end{enumerate}
\end{lemma}

\par\noindent\emph{Proof}. \newline
We try to define by induction on $i\in\omega$ pairs of the form $(D_i,a_i)$ such that $D_i\subseteq C$ is infinite and $a_i\in\omega_2\times\omega$.
We indicate that this attempt is doomed to failure after finitely many steps.

At the stage of $i=0$ we choose any infinite independent $D_0\subseteq C$.
The existence of such a set follows from Ramsey's theorem upon defining $d:[C]^2\rightarrow\{0,1\}$ by $d(x,y)=0$ iff $\{x,y\}\notin E$.
Ramsey's theorem ensures the existence of an infinite monochromatic set $D_0\subseteq C$. The assumption that $G$ has no infinite path implies that $D_0$ must be $0$-monochromatic, that is an independent set.

Now we ask whether there is an element $a\in\omega_2\times\omega$ so that $a$ is connected with $\aleph_0$-many elements from $D_0$.
If the answer is positive then let $a_0$ be the $<_{\rm lex}$-first such element.
If the answer is negative then the process is terminated.

At the stage of $i+1$ we let $D_{i+1}=D_i\cap E(a_i)$.
By the induction hypothesis at the $i$th stage, $D_{i+1}$ is infinite.
Now we ask if there exists some $a\in\omega_2\times\omega$ connected with $\aleph_0$-many elements from $D_{i+1}$ such that $a\neq a_j$ for every $j\leq i$. If yes, let $a_{i+1}$ be the $<_{\rm lex}$-first with this property.
If not, the process is terminated.

Remark that for some $\ell\in\omega$ we will be able to define $D_\ell$ but not $a_\ell$. Otherwise, for every $i\in\omega$ choose an element $d_i\in D_{i+2} - \{d_j:j<i\}$ (here we use the infinitude of each $D_i$) and then $\langle a_n,d_n:n\in\omega\rangle$ forms an infinite path, a contradiction.

Set $D=D_\ell, A=\{a_i:i<\ell\}$.
Define a coloring $c:[D]^2\rightarrow\{0,1\}$ as follows.
Let $c(x,y)=0$ iff there exists $v\notin A$ such that both $\{x,v\},\{y,v\}\in E$.
By another application of Ramsey's theorem there is an infinite $B\subseteq D$ which is monochromatic under $c$.
Observe that $B$ must be $1$-monochromatic, since if $x,y\in B, c(x,y)=0$ then one can produce an infinite path from the elements of $B$ and the elements $v\notin A$ which connect them.
This argument uses the fact that every such $v$ is not in $A$ hence connected with only finitely many elements of $B$.

Now the sets $A,B$ are as required.
First, $A$ is finite and $B$ is infinite.
Second, $(a)$ follows from the choice of the elements of $A$, $(b)$ follows from the fact that $B$ is $c$-monochromatic and $(c)$ from the fact that $B\subseteq D$ and $D$ is independent.

\hfill \qedref{lemfinite}

Equipped with the above lemma, we can proceed to the following.

\begin{claim}
\label{clmclean} Assume that:
\begin{enumerate}
\item [$(a)$] $V\subseteq\omega_2\times\omega$ and ${\rm otp}(V)=\omega^*\cdot\omega_2$.
\item [$(a)$] $G=(V,E)$ is a graph with no infinite path.
\end{enumerate}
Then there exists $W\subseteq V, {\rm otp}(W)=\omega^*\cdot\omega_2$ such that the graph $H=(W,E\cap[W]^2)$ has clean columns.
\end{claim}

\par\noindent\emph{Proof}. \newline
We may assume that all the columns of $V$ are infinite, since ${\rm otp}(V)=\omega^*\cdot\omega_2$ and hence it will remain with the same order type after removing all the finite columns.
We apply Lemma \ref{lemfinite} to every column of $V$, and we get some $U\subseteq V, {\rm otp}(U)=\omega^*\cdot\omega_2$, every nonempty column of $U$ is infinite and edge-free and for each $U(\beta)$ there is a finite set $A(\beta)$ as in the lemma.

Denote the left component $\{\beta: |U(\beta)|=\aleph_0\}$ by $I$, and define $f:I\rightarrow[I]^{<\omega}$ by $f(\beta)=\{\alpha\in I: \exists m\in\omega, (\alpha,m)\in A(\beta)\}$.
Notice that $\beta\notin f(\beta)$ for every $\beta\in I$, since $U(\beta)$ is edge-free and hence no pair of the form $(\beta,m)$ can be an element of $A(\beta)$.
This means that $f$ is a set-mapping.
Further, for every $\beta\in I$ one can see that $f(\beta)$ is a finite set.
This is simply because $A(\beta)$ is finite, due to Lemma \ref{lemfinite}.
By Hajnal's free set theorem there exists $J\subseteq I, |J|=\aleph_2$ such that $J$ is $f$-free, that is $\alpha\notin f(\beta)$ whenever $\alpha,\beta\in J$.

Define $W=\bigcup\{U(\beta):\beta\in J\}$ and observe that ${\rm otp}(W)=\omega^*\cdot\omega_2$. The fact that $H=(W,[W]^2\cap E)$ has clean columns comes from the properties of each $U(\beta)$ as guranteed by the lemma, so we are done.

\hfill \qedref{clmclean}

The ability to clean the columns is helpful in the proof of the main theorem. The proof depends on two additional lemmata.
The second lemma will be postponed after the proof of the main theorem.
For the first lemma let us define the concept of \emph{a replete ordinal}.
Let $G=(V,E)$ be a graph with $V\subseteq\omega_2\times\omega$, and assume that $\psi=\{(\alpha_i,m_i):i\in\omega\}\subseteq V$.
An ordinal $\beta\in\omega_2$ will be called $\psi$-replete iff there exists $n(\beta)\in\omega$ such that for every $k\in[n(\beta),\omega)$ there is $i_k\in\omega$ for which $\{(\alpha_{i_k},m_{i_k}),(\beta,k)\}\in E$.

\begin{lemma}
\label{lemsecond} Suppose that:
\begin{enumerate}
\item [$(a)$] $2^\omega=\omega_1$.
\item [$(b)$] $V\subseteq\omega_2\times\omega$ and ${\rm otp}(V)=\omega^*\cdot\omega_2$.
\item [$(c)$] $H=(V,E)$ is a graph with clean columns.
\item [$(d)$] There is no independent subset of $V$ of type $\omega^*\cdot\omega_2$.
\item [$(e)$] $\psi=\{(\alpha_i,m_i):i\in\omega\}\subseteq V$.
\item [$(f)$] $R$ is an unbounded subset of $\omega_2$ such that every $\beta\in{R}$ is $\psi$-replete.
\end{enumerate}
Then there exists an infinite path in $H$.
\end{lemma}

\par\noindent\emph{Proof}. \newline
For every $\beta\in R$ let $n(\beta)\in\omega$ be as in the definition of repleteness and let $A_\beta\in[\omega]^\omega$ be such that if $i\in A_\beta$ then there is $k\in\omega$ so that $i=i_k$, that is $\{(\alpha_i,m_i),(\beta,k)\}\in E$.
Since $2^\omega=\omega_1$ and $|R|=\aleph_2$ we may assume that $A_\beta=A$ for every $\beta\in R$, where $A$ is some fixed element of $[\omega]^\omega$.
Similarly, we can assume that $n(\beta)\in\omega$ is the same natural number for every $\beta\in R$, and without loss of generality $n(\beta)=0$ for every $\beta\in R$.

Define $c:[R]^2\rightarrow 2$ by $c(\beta,\gamma)=0$ iff there is no edge from $(\beta,k)$ to $(\gamma,\ell)$ whenever $k,\ell\in\omega$.
Put another way, $c(\beta,\gamma)=1$ iff there are $k,\ell\in\omega$ for which $\{(\beta,k),(\gamma,\ell)\}\in E$.
By the Erd\H{o}s-Dushnik-Miller theorem either some $S\in[R]^{\omega_2}$ is $0$-monochromatic or some $\{\beta_m:m\in\omega\}\subseteq R$ is $1$-monochromatic.
In the first case $S\times A$ forms an independent subset of $V$\footnote{We may assume, without loss of generality, that $S\times{A}\subseteq{V}$.} of type $\omega^*\cdot\omega_2$, contradicting $(d)$.
We conclude, therefore, that $\{\beta_m:m\in\omega\}\subseteq R$ is $1$-monochromatic for some infinite subset of $R$.

By induction on $m\in\omega$ we try to choose an element $t_m\in\psi$ such that $\{t_m,(\beta_m,k_m)\}\in E$ and $m<n<\omega\Rightarrow t_m\neq t_n$ and for some $\ell$ we have $\{(\beta_m,k_m),(\beta_{m+1},\ell)\}\in E$.
This is possible since $c(\beta_m,\beta_{m+1})=1$ so we fix $k_m$ and $\ell$ for which $\{(\beta_m,k_m),(\beta_{m+1},\ell)\}\in E$ and then we can choose $t_m$ and $t_{m+1}$ from $\psi$ using the fact that both $\beta_m$ and $\beta_{m+1}$ are $\psi$-replete.\footnote{By a careful choice of $A$ we may assume that $t_m\neq t_{m+1}$.}
Now the sequence $\langle(\beta_m,k_m),(\beta_{m+1},k_{m+1}),t_{m+1}:m\in\omega\rangle$ forms an infinite path in $H$ so the proof is accomplished.

\hfill \qedref{lemsecond}

We can prove now the substantial result which reads as follows:

\begin{theorem}
\label{thmmtgma} Assume $2^\omega=\omega_1, 2^{\omega_1}>\omega_2$ and the generalized Martin's axiom holds.
Then $\omega^*\cdot\omega_2 \rightarrow (\omega^*\cdot\omega_2,\text{infinite path})^2$.
\end{theorem}

\par\noindent\emph{Proof}. \newline
Let $H=(V,E)$ be a graph with no infinite path such that ${\rm otp}(V)=\omega^*\cdot\omega_2$.
We are assuming toward contradiction that there is no independent subset of $V$ of type $\omega^*\cdot\omega_2$.
By Claim \ref{clmclean} we may assume that $H$ has clean columns.
As annotated above, let $I=\{\beta\in\omega_2:|V(\beta)|=\aleph_0\}$.

We define a forcing notion $\mathbb{P}$.
A condition $p\in\mathbb{P}$ is a countable independent subset of $V$.
If $p,q\in\mathbb{P}$ then $p\leq_{\mathbb{P}}q$ iff $p\subseteq q$.
By Lemma \ref{lemmetukan} below, $\mathbb{P}$ satisfies $(c)$ of Theorem \ref{thmgma}.
If $p,q\in \mathbb{P}$ and $p\parallel q$ then $p\cup q\in \mathbb{P}$ and it is a least upper bound by the definition of $\leq_{\mathbb{P}}$.
Similarly, if $(p_j:j\in\omega)$ is $\leq_{\mathbb{P}}$-increasing then $\bigcup_{j\in\omega}p_j\in\mathbb{P}$, being countable and independent, and it is a least upper bound.
Hence $\mathbb{P}$ satisfies the requirements of Theorem \ref{thmgma}.

For every $\beta\in I$ let $D_{\beta}=\{p\in\mathbb{P}:\exists\gamma>\beta, \exists^\infty\ell, (\gamma,\ell)\in p\}$.
We claim that each $D_{\beta}$ is a dense open subset of $\mathbb{P}$.
The fact that $D_{\beta}$ is open follows from the definition, so let us prove density.
Suppose that $p=\{(\alpha_i,m_i):i\in\omega\}\notin D_{\beta}$, but $p\in\mathbb{P}$.
We know, therefore, that $p$ is independent and we observe that assumptions $(a)-(e)$ of Lemma \ref{lemsecond} hold with $p$ here stands for $\psi$ there (note that $(d)$ is our assumption toward contradiction).
Since the conclusion of the lemma fails we see that necessarily assumption $(f)$ of the lemma fails.
That is, the set $R$ of $p$-replete ordinals is bounded in $\omega_2$.
Choose $\gamma\in\omega_2$ such that $\gamma>\beta$ and $\gamma>\sup(R)$.
Namely, $\gamma$ is not $p$-replete and hence $(\gamma,\ell)$ is not connected with an element of $p$ for an infinite set $a$ of $\ell$s.
Define $q=p\cup\{(\gamma,\ell):\ell\in a\}$.
Since $p\leq q\in D_{\beta}$ we see that $D_{\beta}$ is dense.

The collection $\mathcal{D} = \{D_{\beta}: \beta\in I\}$ is of size $\aleph_2$ and $2^{\omega_1}>\omega_2$.
Further, $\alpha<\aleph_2\Rightarrow \alpha^{\aleph_0}<\aleph_2$ since $2^\omega=\omega_1$.
Hence there exists a generic filter $G\subseteq\mathbb{P}$ such that $G\cap D_{\beta}\neq\varnothing$ for every $\beta\in I, n\in\omega$.
Define $W = \bigcup\{p:p\in G\}$, and notice that ${\rm otp}(W)=\omega^*\cdot\omega_2$.
Since every $p\in G$ is independent and $G$ is a directed set, $W$ is independent as well, so we arrived at a contradiction and hence we are done.

\hfill \qedref{thmmtgma}

We accomplish the above proof by the following lemma.

\begin{lemma}
\label{lemmetukan} Let $H=(V,E)$ be a graph with no infinite path, and let $\mathbb{P}$ be the associated forcing whose conditions are countable independent sets.
Assume that $2^\omega=\omega_1$.
Then $\mathbb{P}$ satisfies requirement $(c)$ of Theorem \ref{thmgma}.
\end{lemma}

\par\noindent\emph{Proof}. \newline
Let $\{p_\alpha:\alpha\in\omega_2\}$ be a set of conditions in $\mathbb{P}$.
By induction on $\beta\in\omega_2$ we choose a stationary subset $A_\beta$ of $S^{\omega_2}_{\omega_1}$ such that:
\begin{enumerate}
\item [$(a)$] If $i,j\in A_\beta$ then $p_i\parallel p_j$.
\item [$(b)$] The set $W=S^{\omega_2}_{\omega_1}-\bigcup A_\beta$ is not stationary in $\omega_2$.
\end{enumerate}
We can choose $A_\beta$ by applying Lemma \ref{lemstat} inductively over the set $S^{\omega_2}_{\omega_1}-\bigcup_{\gamma<\beta}A_\gamma$.
Let us describe the construction explicitly.
For $\beta=0$ define $d:[S^{\omega_2}_{\omega_1}]^2\rightarrow{2}$ as follows.
If $i,j\in S^{\omega_2}_{\omega_1}$ then let $d(i,j)=0$ iff $p_i\parallel{p_j}$.
By Lemma \ref{lemstat}, either there is a stationary set $T\subseteq S^{\omega_2}_{\omega_1}$ for which $d''[T]^2=\{0\}$ or an $\omega_1$-sequence of elements of $S^{\omega_2}_{\omega_1}$ for which the range of $d$ on its pairs is constantly one.
The second option is impossible since by an application of Lemma \ref{lemwp} it means that there is an infinite path in $H$, contrary to our assumptions.
Thus, the first option holds.
Set $A_0=T$.

For $\beta>0$, without loss of generality $S^{\omega_2}_{\omega_1}-\bigcup_{\gamma<\beta}A_\gamma$ is stationary.
Apply Lemma \ref{lemstat} in the same way to obtain $A_\beta$.
We emphasize that the choice of the $A_\beta$s depend on $\{p_\alpha:\alpha\in\omega_2\}$ and, moreover, on the enumeration of its elements.
Since we wish to use the argument given in that lemma, for every $\beta\in\omega_2$ there will be an ordinal $\eta_\beta\in\omega_2$ and a fixed sequence $c_\beta\subseteq\eta_\beta$ as described in the proof of Lemma \ref{lemstat}.
Recall that $\eta_\beta<\min(A_\beta)$.
Notice that if $c_\beta=c_\gamma$ then $\eta_\beta=\eta_\gamma$ and $A_\beta\cup A_\gamma$ is linked.\footnote{The adjective \emph{linked} means that every two elements from $A_\beta\cup A_\gamma$ are compatible.}
Hence by taking unions of $A_\beta$s whenever possible we may assume that $\beta\neq\gamma\Rightarrow c_\beta\neq c_\gamma$.

Let $T=\{\min(A_\beta):\beta\in\omega_2\}$.
We claim that $T$ is not stationary.
Suppose not, and define $g:T\rightarrow\omega_2$ by $g(\min(A_\beta))=\eta_\beta$.
This is a regressive function so there is a stationary set $T'\subseteq{T}$ and a fixed $\eta\in\omega_2$ such that $\beta\in T'$ implies $\eta_\beta=\eta$.
Moreover, since each $c_\beta$ is countable we may shrink $T'$ to a stationary set for which all the $c_\beta$s are the same fixed sequence, say $c$.
However, this is impossible since by the choice of the $A_\beta$s if $\beta\neq\gamma$ then $c_\beta\neq c_\gamma$.
We conclude, therefore, that $T$ is not stationary.

Re-enumerate the family of sets $(A_\beta:\beta\in\omega_2)$ by $(S_i:i\in\omega_2)$ in such a way that $i<j\Rightarrow\min(S_i)<\min(S_j)$.
Notice that this implies $i\leq\min(S_i)$ for every $i\in\omega_2$.
Define $f:\omega_2\rightarrow\omega_2$ by letting $f(\alpha)=i$ iff $\alpha\in S_i$.
Let $C$ be a club of $\omega_2$ disjoint from $T\cup{W}$.
If $\alpha\in C\cap S^{\omega_2}_{\omega_1}$ then $f(\alpha)<\alpha$ since $f(\alpha)\leq\alpha$ by the fact that $i\leq\min(S_i)$ and $f(\alpha)=\alpha$ is possible only if $\alpha\in{T}$.
Therefore, $f$ is regressive on $C\cap S^{\omega_2}_{\omega_1}$.
If $\alpha_0,\alpha_1\in C\cap S^{\omega_2}_{\omega_1}$ and $f(\alpha_0)=f(\alpha_1)$ then $\alpha_0$ and $\alpha_1$ belong to the same $S_i$, so $p_{\alpha_0}\parallel p_{\alpha_1}$.
This concludes the proof of the lemma.

\hfill \qedref{lemmetukan}

We can prove, finally, the main result of this paper:

\begin{theorem}
\label{thmmt} Tiltan at $\aleph_1$ is consistent with the positive relation $\omega^*\cdot\omega_1 \rightarrow (\omega^*\cdot\omega_1,\text{infinite path})^2$.
\end{theorem}

\par\noindent\emph{Proof}. \newline
Begin with a universe in which there are enough diamonds so that one can create a tiltan sequence at some stationary set $S\subseteq S^{\omega_2}_{\omega}$, and this instance of tiltan is indestructible upon any $\aleph_1$-complete forcing notion.
The description of this construction and the indestructibility proof appear in \cite{MR1623206}. We may assume that $2^\omega=\omega_1$ along with this construction (actually, this is the natural situation).

We force now the generalized Martin's axiom of Theorem \ref{thmgma} with $2^{\omega_1}>\omega_2$. This is done by an $\aleph_1$-complete forcing notion, hence the tiltan over $S$ is preserved. Observe that this means also tiltan at $\aleph_2$, since $S\subseteq\omega_2$.
By Theorem \ref{thmmtgma} we see that the relation $\omega^*\cdot\omega_2 \rightarrow (\omega^*\cdot\omega_2,\text{infinite path})^2$ holds at this stage.

Now we force with $L\acute{e}vy(\aleph_0,\aleph_1)$.
We claim that the statement of our theorem holds in the generic extension by this collapse.
It is easy to see that the tiltan holds at $\aleph_1$.
Basically, the reason is that unbounded subsets of $\aleph_1$ after the collapse contain unbounded sets of size $\aleph_2$ from the previous stage. The formal proof is elaborated in \cite{MR1623206}.
By a similar argument, the relation $\omega^*\cdot\omega_1 \rightarrow (\omega^*\cdot\omega_1,\text{infinite path})^2$ holds after the collpase.

For this, assume that $G=(V,E)$ is a graph over the set $\omega_1\times\omega$ of order type $\omega^*\cdot\omega_1$.
In the stage just before the collapse there is a graph $H=(W,F)$ such that ${\rm otp}(W)=\omega^*\cdot\omega_2$ and such that $W$ is forced to be a subset of $V$ by the collapse.
We indicate that the copy of $\omega^*$ at each column in the generic extension may issue from a copy of type $\omega_1^*$ before the collapse, in which case we can take a subset of it of order type $\omega^*$.
Now if there is an infinite path in $H$ then it will remain infinite after the collapse as an infinite sequence cannot become finite in the generic extension.
If there is no infinite path in $H$ then there is an independent set of type $\omega^*\cdot\omega_2$, and after the collapse it will be of type $\omega^*\cdot\omega_1$, so we are done.

\hfill \qedref{thmmt}

We conclude with a modest attempt to understand the difference between tiltan and diamond as reflected in our result.
Apparently, tiltan should work in the diamond construction of Baumgartner and Larson from \cite{MR1050558} since one has to add only one edge to each independent set.
Nonetheless, tiltan consistently fails and the reason is that rather than guessing uncountable subsets of $\aleph_1$ we try to guess subsets of $\omega_1\times\omega$.
Moreover, we must make sure that every strongly cofinal subset of $\omega_1\times\omega$ is predicted by strongly cofinal subsets of $\beta\times\omega$ for unboundedly many $\beta\in\omega_1$.

Practically, we translate subsets of $\omega_1\times\omega$ to subsets of $\omega_1$, apply our tiltan sequence, and then send it back to $\omega_1\times\omega$.
But we loose information in this process, and in particular it may happen that a strongly cofinal set will be predicted by narrow sets which will not be strongly cofinal after the translation to $\omega_1\times\omega$.
Larson considered in \cite{MR2279658} a stronger version of tiltan in which strongly cofinal sets are predicted by strongly cofinal tiltan elements.
Unfortunately, such a prediction principle implies the continuum hypothesis (hence diamond)\footnote{Recall that tiltan with the continuum hypothesis imply diamond.} and Larson concluded that this direction will not settle the problem of tiltan and $\omega^*\cdot\omega_1\rightarrow(\omega^*\cdot\omega_1, \text{infinite path})^2$.
On the other hand, this fact motivated our attempt to force the consistency of tiltan with $\omega^*\cdot\omega_1\rightarrow(\omega^*\cdot\omega_1, \text{infinite path})^2$.

It would be interesting to ask a similar question with respect to \emph{superclub}, a prediction principle defined by Primavesi in \cite{primavesi}.
A superclub sequence $\langle S_\alpha:\alpha\in\lim(\omega_1)\rangle$ satisfies $S_\alpha\subseteq\alpha$ for every $\alpha$ and for every $A\in[\omega_1]^{\omega_1}$ one can find $B\subseteq A,|B|=\aleph_1$ such that $\{\delta\in\omega_1:B\cap\delta=S_\delta\}$ is a stationary subset of $\omega_1$.
Superclub holds iff there exists a superclub sequence.
One can see that diamond implies superclub and superclub implies tiltan.
Both implications are irreversible.
In particular, superclub is strictly stronger than tiltan, see \cite{MR3787522} and \cite{1179}.

\begin{question}
\label{qsuperclub} Is it consistent that superclub holds at $\aleph_1$ and $\omega^*\cdot\omega_1 \rightarrow (\omega^*\cdot\omega_1,\text{infinite path})^2$?
\end{question}

We indicate that if one begins with superclub (or even diamond) at $\omega_2$ and collapses $\aleph_1$ then superclub fails at $\aleph_1$ in the generic extension, as follows from \cite{MR3787522}.
In other words, the method of the current paper does not resolve the above problem.

\newpage

\bibliographystyle{alpha}
\bibliography{arlist}

\begin{thebibliography}{EHMR84}

\bibitem[BL90]{MR1050558}
James~E. Baumgartner and Jean~A. Larson.
\newblock A diamond example of an ordinal graph with no infinite paths.
\newblock {\em Ann. Pure Appl. Logic}, 47(1):1--10, 1990.

\bibitem[EHMR84]{MR795592}
Paul Erd{\H{o}}s, Andr{\'a}s Hajnal, Attila M{\'a}t{\'e}, and Richard Rado.
\newblock {\em Combinatorial set theory: partition relations for cardinals},
  volume 106 of {\em Studies in Logic and the Foundations of Mathematics}.
\newblock North-Holland Publishing Co., Amsterdam, 1984.

\bibitem[Gar18]{MR3787522}
Shimon Garti.
\newblock Tiltan.
\newblock {\em C. R. Math. Acad. Sci. Paris}, 356(4):351--359, 2018.

\bibitem[GMS20]{MR4245063}
Shimon Garti, Menachem Magidor, and Saharon Shelah.
\newblock Infinite monochromatic paths and a theorem of
  {E}rdos-{H}ajnal-{R}ado.
\newblock {\em Electron. J. Combin.}, 27(2):Paper No. 2.8, 11, 2020.

\bibitem[GS23]{1179}
Shimon Garti and Saharon Shelah.
\newblock Tiltan and superclub.
\newblock {\em C. R. Math. Acad. Sci. Paris}, 361:853--861, 2023.

\bibitem[Jan17]{comet}
Tove Jansson.
\newblock {\em Comet in Moominland, Translated by Elizabeth Portch}.
\newblock Sort of Books, London, 2017.

\bibitem[Lar90]{MR1050560}
Jean~A. Larson.
\newblock Martin's axiom and ordinal graphs: large independent sets or infinite
  paths.
\newblock {\em Ann. Pure Appl. Logic}, 47(1):31--39, 1990.

\bibitem[Lar06]{MR2279658}
Jean~A. Larson.
\newblock Partition relations on a plain product order type.
\newblock {\em Ann. Pure Appl. Logic}, 144(1-3):117--125, 2006.

\bibitem[Ost76]{MR0438292}
A.~J. Ostaszewski.
\newblock On countably compact, perfectly normal spaces.
\newblock {\em J. London Math. Soc. (2)}, 14(3):505--516, 1976.

\bibitem[Pri11]{primavesi}
Alexander Primavesi.
\newblock Guessing axioms, invariance and {S}uslin trees.
\newblock {\em Ph.D. Thesis, University of East Anglia}, 2011.

\bibitem[She78]{MR0505492}
S.~Shelah.
\newblock A weak generalization of {MA} to higher cardinals.
\newblock {\em Israel J. Math.}, 30(4):297--306, 1978.

\bibitem[She98]{MR1623206}
Saharon Shelah.
\newblock {\em Proper and improper forcing}.
\newblock Perspectives in Mathematical Logic. Springer-Verlag, Berlin, second
  edition, 1998.

\bibitem[Tod21]{MR4217964}
S.~Todorcevic.
\newblock Erdos-{K}akutani phenomena for paths.
\newblock {\em Acta Math. Hungar.}, 163(1):168--173, 2021.

\bibitem[Yei85]{Yeivin}
Israel Yeivin.
\newblock {\em The {H}ebrew language tradition as reflected in the {B}abylonian
  vocalization}.
\newblock Texts and studies XII. The academy of the {H}ebrew language, 1985.

\end{thebibliography}

\end{document}